\documentclass[11pt]{article}
\usepackage{amsfonts,latexsym,amsmath,amscd,geometry}
\usepackage{amssymb}
\usepackage{latexsym}

\newcommand \nc{\newcommand}
\newtheorem{theorem}{Theorem}[section]
\newtheorem{lemma}[theorem]{Lemma}

\newtheorem{remark}[theorem]{Remark}

\nc{\ba}{\begin{array}}\nc{\ea}{\end{array}}
\nc{\be}{\begin{eqnarray}}\nc{\ee}{\end{eqnarray}}
\nc{\beq}{\begin{equation}}\nc{\eeq}{\end{equation}}
\nc{\bex}{\begin{eqnarray*}}\nc{\eex}{\end{eqnarray*}}
\nc{\btm}{\begin{theorem}} \nc{\etm}{\end{theorem}}
\nc{\blm}{\begin{lemma}} \nc{\elm}{\end{lemma}}
\nc{\R}{\mathbb{R}}  \nc{\ld}{\lambda}
\nc{\va}{\varphi}
\nc{\ve}{\varepsilon}
\def\x{\mathbf{x}}\def\e{\mathbf{e}}

\def\pa{\partial}
\def\pf{\noindent{\bf Proof.\quad}}\def\endpf{\hfill$\Box$}

\def\D{\widetilde{D}}\def\del{\widetilde{\Delta}}

\begin{document}
\title{Nonuniqueness of nematic liquid crystal flows in dimension three}

\author{Huajun Gong\footnote{College of Mathematics and Statistics, Shenzhen University, Shenzhen 518060, China. huajun84@hotmail.com}, \quad Tao Huang\footnote{ NYU-ECNU Institute of Mathematical Sciences at NYU Shanghai,
3663 Zhongshan Road North, Shanghai, 200062, China. th79@nyu.edu} \quad and\quad Jinkai Li\footnote{Department of Mathematics, The Chinese University of Hong Kong, Shatin, N.T., Hong Kong. jklimath@gmail.com}\\
}
\date{}
\maketitle

\begin{abstract}
For suitable initial and boundary data, we construct infinitely many weak solutions to the nematic liquid crystal flows in dimension three. These solutions are in the axisymmetric class with bounded energy and ``backward bubbling'' at a large time.
\end{abstract}

\section{Introduction}
For any smooth domain $\Omega\subset R^3$, we consider the following  simplified model of nematic liquid crystal flows
\begin{equation}\label{liquidcrystal}
\begin{cases}
u_t+u\cdot\nabla u-\mu\Delta u+\nabla P=-\lambda\nabla\cdot\big(\nabla d\odot\nabla d-\frac12{|\nabla d|^2}\mathbb{I}_3\big),\\
\nabla \cdot u=0, \\
d_t+u\cdot\nabla d=\gamma\big(\Delta d+|\nabla d|^2d\big),
\end{cases}
\end{equation}
where $u(\x,t):\Omega\times(0,T)\rightarrow\R^3$ is the velocity field of the underlying incompressible fluid,
$d(\x,t):\Omega\times(0,T)\rightarrow \mathbb S^2:=\big\{\x\in\R^3:  |\x|=1\big\}$ represents the (averaged) orientation field of nematic liquid crystal molecules,
$P(\x,t):\Omega\times(0,T)\rightarrow\R$ is the pressure function, $\x=(x,y,z)\in\Omega$,
$\nabla\cdot$ denotes the divergence operator on $\mathbb R^3$, $\nabla
d\odot\nabla d =\big(\langle\frac{\partial d}{\partial
\x_i},\frac{\partial d}{\partial \x_j}\rangle\big)_{1\leq i,j\leq 3}\in \R^{3\times 3}$
represents the stress tensor induced by the orientation field $d$, and
$\mathbb{I}_3=\big(\delta_{ij}\big)_{1\le i, j\le 3}\in\R^{3\times 3}$ is the identity matrix of order $3$. The parameters $\mu$, $\lambda$ and $\gamma$ are positive constants representing the fluid viscosity,
the competition between kinetic energy and potential energy, and the macroscopic elastic relaxation time for the molecular orientation field respectively. For simplicity, we only consider
$$
\mu=\lambda=\gamma=1.
$$
The system \eqref{liquidcrystal} has first been proposed by Lin \cite{lin} as a simplified version of the general Ericksen-Leslie system modeling the hydrodynamic flow of nematic liquid crystal materials proposed by Ericksen \cite{ericksen} and Leslie \cite{leslie} between 1958 and 1968. The system (\ref{liquidcrystal}) is a macroscopic continuum description of the time evolution of the liquid crystal material under the influence of both the fluid field and the macroscopic description of the microscopic
orientation configurations of rod-like liquid crystals molecules.
The interested readers can refer to \cite{ericksen}, \cite{leslie}, \cite{lin}, and Lin-Liu \cite{lin-liu} for more details.
In this paper, we will investigate the system \eqref{liquidcrystal} with initial data
$$
(u,d)\big|_{t=0}=(u_0,d_0)
$$
and boundary data that will be specific later,
where $(u_0(\x),d_0(\x)):\Omega\rightarrow \R^3\times\mathbb S^2$ satisfies $\mbox{div}~ u_0=0$ and $|d_0|=1$.

Mathematically, the system \eqref{liquidcrystal} is a strong coupling  between the incompressible Naiver-Stokes equation for the flow field
and the (transport) heat flow of harmonic maps for the
orientation field of the liquid crystal molecules, which has attracted a lot of interests and generated
many interesting results recently. In dimension two, Lin-Lin-Wang \cite{lin-lin-wang} have proved the existence of global Leray-Hopf type weak solutions to initial and boundary value problem of \eqref{liquidcrystal} with finitely many possible singular times
(see \cite{hong} for $\Omega=\R^2$, \cite{HX,WMWWD,LX} for more general systems, and \cite{dong-lei,LLZ,JKL2D,xu-zhang} for some other related works). Lin--Wang \cite{lin-wang}, Wang--Wang--Zhang \cite{WMWWDZZF}, and Li--Titi--Xin \cite{LTX} have also proved the uniqueness for such weak solutions. In dimension three, Lin-Wang \cite{lin-wang3} have proved the existence of global weak solutions under the assumption $d_0(x)\in \mathbb S^2_+$ for a.e. $x\in\Omega$ by developing some new compactness arguments. Here $
\mathbb S^2_+$ is the upper hemisphere. Huang \cite{huang16} has shown the weak solution is regular and unique in the scaling invariant Leray spaces.
Recently, in \cite{hllw16}, two nontrivial examples of finite time singularities in dimension three have been constructed. The first example is built with the help of axisymmetric solutions to \eqref{liquidcrystal} without swirl. In the second example, the initial data of the approximate harmonic maps was constructed with small energy but large topology. With help of the energy inequality, the local smooth solutions have been proved to have finite time singularities by an $\epsilon$-apriori estimate on approximate harmonic maps. However, it still remains a very challenging open problem to establish the existence of global Leray-Hopf type weak solutions and partial regularity of suitable weak solutions to \eqref{liquidcrystal} in dimension three. It should be mentioned that for suitably regular initial data, local existence and uniqueness of more regular solutions than the weak ones can be established for more general systems than (\ref{liquidcrystal}), see, e.g., \cite{HLX,HIEBER1,HIEBER2,GLXNA,GLXZ,JKLM,MGL}. More results and references can be found in the survey paper by Lin-Wang \cite{lin-wangs}.

In order to investigate the nonuniqueness of system \eqref{liquidcrystal}, it is helpful to mention a few of related results for the heat flow of harmonic maps.
In dimension two, a unique global weak solutions with finitely many singularities has been constructed by Struwe \cite{struwe} and Chang \cite{chang}. Freire \cite{freire} proved that the weak solution is uniqueness if the Dirichlet energy is monotone decreasing with respect to $t$ (see also \cite{luwang} and \cite{linlz} for alternative proofs). The examples of non-unique weak solutions, whose Dirichlet energies are not monotone decreasing, have been constructed by Bertsch-Dal Passo-Hout \cite{bph02} and \cite{topping02}. In higher dimensions, the existence of global weak solutions with partial regularity has been established by Chen-Struwe \cite {chen-struwe} and Chen-Lin \cite{chen-lin}. Examples of non-unique weak solutions have been constructed by Coron \cite{coron} and Bethuel-Coron-Ghidaglia-Soyeur \cite{BBC}. Recently, Huang-Wang \cite{hw16} have established several new results on uniqueness in higher dimensions. More references on the heat flow of harmonic maps can be found in the book by Lin-Wang \cite{lin-wangbk}.

Inspired by all the previous results, it is quit an interesting question to investigate the nonuniqueness of (nontrivial) weak solutions to the nematic liquid crystal flows \eqref{liquidcrystal}. More precisely, we would like to construct more than one weak solution with special choice of initial and boundary values with axisymmetric structure.

\begin{theorem}\label{mainth1}
Let $B^2_1$ be the unit disk in dimension two and $\Omega=B_1^2\times [0,1]$ be the round cylinder. For the initial data $(u_0, d_0)\in C^{\infty}(\bar\Omega, \R^3\times\mathbb S^2)$ with
\beq\label{initialu0}
u_0({\bf x})=(x,y,-2z),
\eeq
and
\beq\label{initiald0}\displaystyle d_0({\bf x})=\Big(\frac{x}{\sqrt{x^2+y^2}}\sin\varphi_0\big(\sqrt{x^2+y^2}\big),
\frac{y}{\sqrt{x^2+y^2}}\sin\varphi_0\big(\sqrt{x^2+y^2}\big), \cos\varphi_0\big(\sqrt{x^2+y^2}\big)\Big),
\eeq
for some $\va_0\in C^{\infty}[0,1]$ and $\va_0^{(2k)}(0)=0$, $k\in\mathbb N$, the system \eqref{liquidcrystal} has more than one weak solution satisfying \eqref{liquidcrystal} in the sense of distribution and initial and boundary conditions \eqref{initialu0}-\eqref{initiald0} in the sense of trace.
\end{theorem}

 We also want to mention that, in dimension two \cite{dong-lei} have constructed a family of exact axisymmetric solutions to \eqref{liquidcrystal}, and \cite{vdh01} and \cite{hv01} have also considered other simplified models of liquid crystal flows with cylindrical symmetry and showed some finite time blowup and nonuniqueness results. Recently, \cite{cky16} have also constructed solutions to system \eqref{liquidcrystal} in dimension three that are twisted and periodic along $z$-axis.

The idea of our proof is mainly inspired by \cite{bph02} and \cite{topping02}. The main difficulty for our case is that there is no global existence for liquid crystal flows \eqref{liquidcrystal} in dimension three with general initial and boundary conditions (except for those in \cite{lin-wang3} in our best knowledge). However, very fortunately, after the first blowup at $t_1$ in our construction, the value of the solution at singular point will jump from $0$ to $\pi$, which will make the solution always topologically small after the singular time. Thus by the recent existence result in \cite{hllw16}, we will have one global weak solution to \eqref{liquidcrystal} with single singularity (Theorem \ref{thm3.2}).

The solution will become no bigger than $\pi$ at some time $t_2$ while the value of the solution at origin will be kept as $\pi$ (Theorem \ref{thm4.1}). But to keep it as $\pi$ at origin, some extra energy seems be held until $t_2$, which may be released at any time after $t_2$ if the value of the solution at origin drops to zero again. This will produce an possible upwards jump of the energy at any time after $t_2$. Therefore, we may have infinitely many weak solutions to the original system \eqref{liquidcrystal} (Theorem \ref{thm5.1}).

The paper is organized as follows. In Section 2, we will review a simplified form of the axisymmetric liquid crystal flows that has been first considered in \cite{hllw16}. In Section 3, we will construct the first global weak solution with  finite time blowup and energy drop. In Section 4, by construction of a suitable supersolution, we will construct a crucial time $t_2$ in finding another solution. In Section 5, we will construct a upward energy jump beyond $t_2$ and prove our main results.

\section {Axisymmetric solutions without swirls}
\setcounter{equation}{0}
\setcounter{theorem}{0}

To construct weak solutions, we first introduce the axisymmetric solutions to (\ref{liquidcrystal}) without swirls, i.e.,
$$\begin{cases} u(r,\theta, z,t)=u^{r}(r,z,t)\e^r+u^{3}(r,z,t)\e^3,\\
d(r,\theta, z,t)=\sin\va(r,z,t) \e^r+\cos\va(r,z,t) \e^3,\\
P(r,\theta,z,t)=P(r,z,t).
\end{cases}
$$
A domain $\Omega\subset\mathbb R^3$ is axisymmetric if it is invariant under any rotation map. Then $(u^r, u^3, \varphi, P)$ solves (cf. \cite{hllw16})
\begin{equation}\label{liquidcrystal1}\displaystyle
\begin{cases}
\displaystyle\frac{\D u^r}{Dt}-\mu\del u^r+\frac{u^r}{r^2}+P_r
=-\big(\del\va-\frac{\sin(2\va)}{2r^2}\big)\va_r,\\
\ \ \ \ \ \ \ \ \ \displaystyle\frac{\D u^3}{Dt}-\del u^3+P_z
=-\big(\del\va-\frac{\sin(2\va)}{2r^2}\big)\va_z,\\
\qquad\qquad\ \displaystyle\frac{1}{r}(ru^r)_r+(u^3)_z=0,\\
\qquad\qquad\qquad\ \ \displaystyle\frac{\D \va}{Dt}-\del\va
=-\frac{\sin(2\va)}{2r^2},
\end{cases}
\end{equation}
where
\bex
\displaystyle\frac{\D}{Dt}:=\pa_t+u^r\pa_r+u^3\pa_z,\qquad \del:=\pa^2_r+\frac{1}{r}\pa_r+\pa^2_z,
\eex
are the material derivatives and Laplace operator in cylindrical coordinates respectively.
According to \cite{hllw16}, to construct a special nontrivial solution of Navier-Stokes equations, we will only consider the round cylinder $\Omega=B_1^2\times [0,1]$
and an axisymmetric solution $(u,P, d)$ without swirl in the following form
$$\begin{cases}
u(r,\theta, z,t):=v(r,t)\e^r+w(z,t)\e^3,\\
d(r,\theta, z, t):=\sin\va(r,t) \e^r+\cos\va(r,t) \e^3,\\
P(r,\theta, z,t):=Q(r,t)+R(z,t).
\end{cases}
$$
Then \eqref{liquidcrystal1} becomes
\begin{equation}\label{liquidcrystal2}
\begin{cases}
\displaystyle v_t+vv_r-\big( v_{rr}+\frac{v_r}{r}-\frac{1}{r^2}v\big)+Q_r
=\big(\va_{rr}+\frac{\va_r}{r}-\frac{\sin(2\va)}{2r^2}\big)\va_r, \  r\in [0,1], \\
\qquad\qquad\quad \ w_t+ww_z- w_{zz}+R_z
=0, \ \qquad\qquad \ \ \ \ \ \ \ \ \ \ \ \ \ \ \ \ \ \ z\in [0,1], \\
\qquad\qquad\qquad\qquad\ \ \ \ \ \ \displaystyle\frac{1}{r}(rv)_r+w_z=0,
\ \ \ \ \ \ \ \ \ \ \ \ \ \ \ \    (r,z)\in [0,1]\times [0,1],\\
\qquad\qquad\ \ \ \ \ \ \ \ \ \ \ \ \ \ \ \ \ \ \ \ \ \  \displaystyle\va_t+v\va_r=\va_{rr}+\frac{\va_r}{r}-\frac{\sin(2\va)}{2r^2},
 \qquad\ r\in [0,1].
\end{cases}
\end{equation}
If we only consider the special initial conditions of $(u_0, d_0)$ given by \eqref{initialu0} and \eqref{initiald0},  the initial condition of $(v,w,\varphi)$ is given by
\begin{equation}\label{sbinitial}
\begin{cases}
v|_{t=0}=r,\ \ \ \ \ \ \ 0\le r\le 1,\\
w|_{t=0}=-2z, \ \ 0\le z\le 1,\\
\va|_{t=0}=\va_0(r),\ 0\le r\le 1,
\end{cases}
\end{equation}
for some $\va_0\in C^\infty([0,1])$, with $\va_0(0)=0$.
The boundary conditions are
\beq\label{sboundary}
\begin{split}
\begin{cases} v(0,t)=0\\
v(1,t)=1,
\end{cases} & \quad \begin{cases} w(0,t)=0\\ w(1,t)=-2,\end{cases}
\quad \begin{cases} \va(0,t)=\chi(t)\\ \va(1,t)=\va_0(1).\end{cases}
\end{split}
\eeq
Here the $\chi(0)=0$, $\chi(t)=k(t)\pi$ and $k(t)$ equals to $0$ or $1$ depending on time.

According to the Lemma 2.3 and Lemma 2.4 in Huang-Lin-Liu-Wang \cite{hllw16} and combining the the initial and boundary condition (\ref{sbinitial})-\eqref{sboundary} with divergence free condition, we conclude for any $(r,z,t)\in [0,1]\times [0,1]\times [0,T)$
\beq
v(r,t)=r,\quad w(z,t)=-2z
\eeq
and
 \begin{equation}\label{pressure}
\displaystyle R(z,t)=-2z^2+c_1(t),\quad
Q(r,t)
=-\int_0^r\big(\va_{rr}+\frac{\va_r}{r}-\frac{\sin(2\va)}{2r^2}\big)\va_r\,dr-\frac{r^2}2+c_2(t),
\end{equation}
for some $c_1, c_2\in C^\infty([0,T))$.
Therefore, there only equation left is the following drift heat flow of axisymmetric harmonic flows in dimension two
\begin{equation}\label{drift-hm}
\begin{cases}
\displaystyle\va_t+r\va_r=\va_{rr}+\frac{\va_r}{r}-\frac{\sin(2\va)}{2r^2}, \ 0<r<1, \\
\varphi(r,0)=\varphi_0(r), \ 0<r<1,\\
\varphi(0,t)=\chi(t), \ \ \varphi(1,t)=\varphi_0(1), \ t>0.
\end{cases}
\end{equation}
Now, it is easy to see that Theorem \ref{mainth1} will be indicated by the following result.

\begin{theorem}\label{mainth2}
There exists a $\va_0\in C^{\infty}[0,1]$ with $\va_0^{(2k)}(0)=0$, $k\in \mathbb N$, such that the
drift heat flow of axisymmetric harmonic flows \eqref{drift-hm} has more than one weak solution.
\end{theorem}
The local existence of unique smooth solution to \eqref{drift-hm} has been proved by \cite{hllw16} with $\chi(t)\equiv0$. We will extend this local solution in different ways after the first blowup by choosing different $\chi(t)$ and the exact proof of Theorem \ref{mainth2} will be given in Section 5.

We also want to mention that if $\varphi$ solves  (\ref{drift-hm}),
then
$$d(x,y,t):=\Big(\sin\va\big(\sqrt{x^2+y^2},t\big)\frac{x}{\sqrt{x^2+y^2}},
\sin\va\big(\sqrt{x^2+y^2},t\big)\frac{y}{\sqrt{x^2+y^2}},
\cos\va\big(\sqrt{x^2+y^2},t\big)\Big),$$
for $(x,y)\in B_1^2$ and $t\ge 0$, solves
\beq\label{eqd2}
d_t+(x,y)\cdot\nabla d=\Delta d+|\nabla d|^2 d,\ (x,y)\in B_1^2,\ t>0,
\eeq
with the initial condition
\beq\label{d2initial}
d\big|_{t=0}=\Big(\sin\va_0\big(\sqrt{x^2+y^2}\big)\frac{x}{\sqrt{x^2+y^2}},
\sin\va_0\big(\sqrt{x^2+y^2}\big)\frac{y}{\sqrt{x^2+y^2}}, \cos\va_0\big(\sqrt{x^2+y^2}\big)\Big),
\eeq
and the boundary condition
\beq\label{d2boundary}
d\big|_{\partial B_1^2}=\big(\sin \va_0(1) x, \sin\va_0(1)y, \cos \va_0(1) \big).
\eeq
By Remark 5.3 in \cite{hllw16}, in general the liquid crystal equations \eqref{liquidcrystal1} with initial and boundary conditions (\ref{sbinitial})-\eqref{sboundary} may not satisfy the energy inequality. But if we only consider the equation \eqref{eqd2} with initial and boundary conditions \eqref{d2initial} and \eqref{d2boundary}, some special form of energy inequality will be valid.

\begin{lemma}\label{lemma2.2}
For any $0<T<+\infty$ and any regular solutions to the equation \eqref{eqd2} with initial and boundary conditions \eqref{d2initial} and \eqref{d2boundary}, it holds
\beq
\int_{B_1^2}|\nabla d|^2(\cdot, T)\,d\x
\leq e^{T}\int_{B_1^2}|\nabla d_0|^2\,d\x
\eeq
and
\beq
\int_{0}^{T}\int_{B_1^2}|d_t|^2\,d\x dt
\leq C(T).
\eeq
\end{lemma}

\begin{remark}\label{rmk2.3}
By the lemma, it is easy to see
\beq
\int_{0}^{T}\int_{0}^1|\va_t|^2r\,dr dt
+\int_{0}^1\left(|\va_r|^2+\frac{\sin^2 \va}{r^2}\right)r\,dr
\leq C(T).
\eeq

\end{remark}

\pf Multiplying both sides of equation \eqref{eqd2} by $d_t$, integrating in $[0,1]$ and using integration by parts and Young's inequality, one should have
\beq\notag
\begin{split}
\frac12\frac{d}{dt}\int_{B_1^2}|\nabla d|^2\,d\x+\int_{B_1^2}|d_t|^2\,d\x=-\int_{B_1^2}(x,y)\cdot\nabla d d_t\,d\x
\leq \frac12\int_{B_1^2}|d_t|^2\,d\x+\frac12\int_{B_1^2}|\nabla d|^2\,d\x.
\end{split}
\eeq
Thus, the energy estimates are implied by Gronwall's inequality.
\endpf

\section{Construction of first solution with single blowup}
\setcounter{equation}{0}
\setcounter{theorem}{0}

We consider the initial data $\va_0(r)\in C^{\infty}[0,1]$ satisfying the following properties:
\beq\label{asspva0}
\left\{\begin{array}{ll}
\vspace{2mm}
\displaystyle\va_0^{(2k)}(0)=0,\ k\in\mathbb N,\quad \va_0\left(\frac12\right)=\alpha\in (\pi,2\pi),\quad
\va_0\left(1\right)=\beta\in (0,\pi),\quad\\
\displaystyle\va_0\mbox{ is increasing in } (0,\frac12)\mbox{ and decreasing in }(\frac12,1).
\end{array}
\right.
\eeq
By \cite{hllw16}, the local existence of smooth axisymmetric solution to \eqref{drift-hm} has been constructed for smooth initial data $\va_0$. In order to construct finite time blowup, we need the following useful comparison principle for the weak sub(sup)-solutions inspired by \cite{bph02}. Denote
$$
M_t:=(0,1)\times(0,t),\quad \partial M_t=[0,1]\times \{0\}\cup \{0,1\}\times[0,1).
$$

\begin{lemma}\label{complm}
Suppose $f,g\in C(M_t)\cap H^1_{\footnotesize\mbox{loc}}(M_t)$ satisfy the following inequalities
\beq
\iint_{M_t}(f_t+rf_r)\xi r\,drdt\leq -\iint_{M_t}\left(f_r \xi_r+\frac{\sin 2f}{2r^2}\xi\right)r\,drdt,
\eeq
\beq
\iint_{M_t}(g_t+rg_r)\xi r\,drdt\geq -\iint_{M_t}\left(g_r \xi_r+\frac{\sin 2g}{2r^2}\xi\right)r\,drdt,
\eeq
for any $0\leq \xi\in C^{\infty}_0(M_t)$. If $f\leq g$ on $\partial M_t$, then it holds $f\leq g$ in $\overline{M_t}$.
\end{lemma}
The proof of the lemma is the same as in \cite{bph02}, so we omit it here. As stated in Appendix of \cite{bph02}, the same conclusion holds when $f,g$ are allowed to be discontinuous at finitely may points $\{(0,\bar t_i)\}$ for $0<t_i\leq t$.

Our main result of this section can be stated as follows.

\begin{theorem}\label{thm3.2}
There exists a function $\va_0$ satisfying the assumption \eqref{asspva0} such that the system \eqref{drift-hm} has a weak solution $\tilde \va$ which is smooth in $[0,1]\times[0,+\infty)$ with exception of one point $(0,t_1)$ for some $0<t_1<+\infty$. It also holds
\beq
\lim\limits_{r\rightarrow 0} \tilde\va(r,t_1)=\pi,\quad
\tilde\va(0,t)=\left\{
\begin{array}{ll}
0,\quad \mbox{ if  }\ 0\leq t<t_1,\\
\pi,\quad \mbox{ if  }\ t>t_1,
\end{array}
\right.
\eeq
and
\beq\label{engdrop1}
\lim\limits_{t\uparrow t_1} E(\tilde\va(t))\geq 4+E(\tilde \va(t_1)),
\eeq
where
$$
E(\tilde\va(t))=\int_{0}^1\left(|\tilde\va_r|^2+\frac{\sin^2 \tilde\va}{r^2}\right)r\,dr.
$$
\end{theorem}

\pf We divide our proof into several steps.

\noindent{\bf Step 1.} We first construct the local smooth solution with finite time blowup with a help of the subsolution as in \cite{hllw16}. Let
\beq\label{subsol}
\Phi(r,t)=2\arctan\left(\frac{r}{\lambda(t)e^t}\right)
+2\arctan\left(\frac{r^{1+\ve}}{\mu e^t}\right),
\eeq
where $\lambda(t)$ is a solution of the ordinary differential equation
$$
\lambda'=-\delta e^{-2t} \lambda^\ve, \mbox{ for } t>0,\quad \lambda(0)=\lambda_0>0.
$$
For $0<2\lambda_0^{1-\ve}<\delta(1-\ve)$, denote
$$
T_{\lambda}:=\frac{1}{2}\ln \left(\frac{\delta(1-\ve)}{\delta(1-\ve)-2\lambda_0^{1-\ve}}\right)>0.
$$
be the first time when $\lambda(t)=0$. Then there exist suitable $\mu>0$, $\ve\in(0,1)$ and $\delta>0$ (cf. \cite{hllw16}) such that
$\Phi\in C^{\infty}([0,1]\times[0,T_\lambda))$ satisfies
\beq\label{bdyPhi}
\lim\limits_{r\rightarrow 0}\Phi(r,t)=0,\mbox{ for } t\in[0,T_\lambda),
\quad
\lim\limits_{r\rightarrow 0}\Phi(r,T_\lambda)=\pi,
\eeq
and
\beq\label{eqnPhi}
\Phi_{rr}+\frac{\Phi_r}{r}-\frac{\sin 2\Phi}{2r^2}-r\Phi_r-\Phi_t\geq 0
\eeq
in $(0,1)\times(0,T_\lambda)$.

Let $\underline \va_0(r)$ be a function satisfies the assumption \eqref{asspva0}. Then the system \eqref{drift-hm} with initial data $\underline \va_0$ and $\chi(t)=0$ has a unique classic solution
$\underline \va(r,t)$ in $[0,1]\times[0,\underline t_0)$ for some $\underline t_0>0$. The fact $\pi<\underline \va_0\left(\frac12\right)<2\pi$ implies that there exist $\underline t^*<\underline t_0$ and $\gamma>\pi$ such that
$$
\underline \va(\frac{1}{2}, t)\geq \gamma, \quad  \mbox{ for }0\leq t\leq\underline t^*.
$$
Choosing suitable $\delta$ and $\lambda_0$ such that $T_\lambda<\underline t^*$ and $\Phi(r,t)$ is a subsolution to the following problem
\begin{equation}\label{eqnsubsl}
\begin{cases}
\displaystyle f_t+rf_r=f_{rr}+\frac{f_r}{r}-\frac{\sin(2f)}{2r^2},\quad \ 0<r<\frac12,\ t>0,\\
f(r,0)\geq\Phi(r,0),\quad \ 0<r<\frac12,\\
f(0,t)=0,\quad \ \ \ f\left(\frac12,t\right)=\gamma, \quad\ t>0.
\end{cases}
\end{equation}
Now, we can choose another initial data $\va_0(r)$ satisfying \eqref{asspva0}  and
$$
\va_0(r)\geq \underline \va_0(r),\ \ \mbox{ for } r\in [0,1],\quad\mbox{and}\quad \va_0(r)\geq\Phi(r,0),\ \ \mbox{for } r\in \left[0,\frac12\right].
$$
By the comparison principle Lemma \ref{complm} (see also \cite{hllw16}) and \eqref{bdyPhi}-\eqref{eqnPhi}, there exists a time $0<t_1\leq T_\lambda\leq \underline t^*$ such that the system \eqref{drift-hm} with initial data $\va_0$ and $\chi(t)=0$ has a unique classic solution
$ \va(r,t)$ in $[0,1]\times[0,t_1)$, which blows up at $(r,t)=(0,t_1)$ and
\beq\label{nuq2.7}
\limsup\limits_{(r,t)\rightarrow (0,t_1)} \va(r,t)>0.
\eeq

\vspace{2mm}
\noindent{\bf Step 2.} We now proceed to find the limit of $\va(r,t)$ as $(r,t)\rightarrow (0,t_1)$. Direct calculation implies
$$
\bar\phi(r,t)=2\arctan \left(\frac{r}{\sigma e^t}\right)+\pi,
$$
is a solution to the equation
$$
\bar\phi_t+r\bar\phi_r-\bar\phi_{rr}-\frac{\bar\phi_r}{r}+\frac{\sin 2\bar\phi}{2r^2}=0.
$$
Choose the constant $\sigma$ large enough such that $\va_0(r)\leq \bar\phi(r,0)$ for any $r\in[0,1]$. Then by the comparison principle Lemma \ref{complm} (see also \cite{hllw16}), we obtain
$$
\va(r,t)\leq \bar\phi(r,t),\quad\mbox{for }(r,t)\in [0,1]\times[0,t_1),
$$
which implies
$$
\limsup\limits_{(r,t)\rightarrow (0,t_1)} \va(r,t)
\leq \lim\limits_{(r,t)\rightarrow (0,t_1)} \bar\phi(r,t)=\pi.
$$
Define
$$
r^*(t)=\sup\left\{r>0~\big|~0\leq \va(r',t)<\pi,\ \mbox{for all }0\leq r'\leq r\right\}.
$$
If we assume $r^*(t)\geq \delta_0>0$ for all $0<t<t_1$, then there exists a $0<\delta_0'<\delta_0$ such that
$$
\va(r,t)\leq \pi-c_0,\quad \mbox{in }[0,\delta_0']\times[0,t_1)
$$
for some constant $c_0>0$. By \cite{hllw16}, we can construct a smooth supersolution $\tilde\phi(r,t)$ to the system with $\tilde\phi(0,t)=0$ for any $0\leq t\leq t_1$, which is contradict to \eqref{nuq2.7} and implies
\beq\label{nuq2.9}
\liminf\limits_{t\rightarrow t_1} r^*(t)=0,
\eeq
Therefore, we obtain
\beq\label{nuq2.6}
\limsup\limits_{(r,t)\rightarrow (0,t_1)} \va(r,t)=\pi.
\eeq
By the definition of $r^*$ and \eqref{nuq2.9}, in order to argue
\beq\label{nuq2.10}
\lim\limits_{r\rightarrow 0}\va(r,t_1)=\pi,
\eeq
we only need to show that
$$
h(r,t)\geq \pi,\quad \mbox{for any }r\in[r^*,\frac{1}{2}],\ t\in (0,t_1).
$$
Indeed, denote $\omega(r,t)=h(r,t)-\pi$. It is easy to see that $\omega$ is continuous in $\displaystyle\left[0,\frac12\right]\times[0,t_1)$ and is a classic solution of
$$\omega_t+r\omega_r=\omega_{rr}+\frac{\omega_r}{r}-\frac{\sin 2\omega}{2r^2},\quad \mbox{in }\left(0,\frac{1}{2}\right)\times(0,t_1).$$
The initial and boundary data satisfy
$$
\omega_0(r)<0,\quad \mbox{if }r\in [0,r^*),\quad \omega_0(r)>0,\quad \mbox{if }r\in (r^*,\frac12),
$$
and
$$
\omega(0,t)<0<\omega\left(\frac12,t\right), \quad \mbox{for }t\in[0,t_1).
$$
By the maximum principle and the Jordan curve theorem as in \cite{natano82} (Lemma 2.5 and 2.6), we conclude that $\omega(r)$ will keep positive in $[r^*,\frac12]$, which completes the proof of \eqref{nuq2.10}.

\vspace{2mm}
\noindent{\bf Step 3.} To show the energy drop inequality \eqref{engdrop1}, we first denote
$$
E_{r_1,r_2}(f):=\int_{r_1}^{r_2}\left((f')^2+\frac{\sin^2 f}{r^2}\right)\,rdr,
$$
and
$$
W_r:=\big\{f\in H^1_r(0,r)~\big|~ E_{0,r}(f)<\infty,\ f(0)=0,\ f(r)=h(r,t)\big\},
$$
where $X_r$ denotes Sobolev spaces $X$ with weight $r$. By the result in \cite{vdh01}, we obtain if $|\va(r_2,t)|<\pi$,
$$
E_{0,r_2}(\va(t))\geq \min\limits_{f\in W_{r_2}} E_{0,r_2}(f)
\geq 2(1-\cos \va(r_2,t)).
$$
Thus if we choose $\tilde r\in (0,1)$ small enough, we have
\beq\notag
\begin{split}
&\lim\limits_{t\uparrow t_1} E(\va(t))\\
=&\lim\limits_{t\uparrow t_1}\big(E_{0,\tilde r}(\va(t))+E_{\tilde r,1}(\va(t))\big)\\
\geq& \lim\limits_{t\uparrow t_1}\big(2(1-\cos \va(\tilde r, t))+E_{\tilde r,1}(\va(t))\big)\\
=&2(1-\cos \va(\tilde r,t_1))+E_{\tilde r,1}(\va(t_1)).
\end{split}
\eeq
By \eqref{nuq2.10}, we conclude that
$$
\lim\limits_{t\uparrow t_1} E(\va(t))\geq E(\va(t_1))+4,
$$
which is exact the energy drop inequality \eqref{engdrop1}.

\vspace{2mm}
\noindent{\bf Step 4.} To extend the solution $\va(r,t)$ beyond $t_1$, we consider the system
\begin{equation}\label{nuq2.8}
\begin{cases}
\displaystyle \va_t+r\va_r=\va_{rr}+\frac{\va_r}{r}-\frac{\sin(2\va)}{2r^2},\quad \ 0<r<1,\ t>0,\\
\va(r,0)=\va(r,t_1)-\pi,\quad \ 0\leq r\leq 1,\\
\va(0,t)=0,\quad \ \ \ \va\left(1,t\right)=\beta-\pi, \quad\ t>0.
\end{cases}
\end{equation}
It is easy to see that $\va(r,t_1)-\pi\in C^0[0,1]$ and $|\va(r,t_1)-\pi|\leq \pi-c_1$ with some constant $c_1>0$. Thus by the global existence result in \cite{hllw16}, the system \eqref{nuq2.8} has a classic solution $\va_1(r,t)$ for $(r,t)\in (0,1)\times(t_1,+\infty)$. We can complete the proof of Theorem \ref{thm3.2} by setting
\beq
\tilde\va(r,t):=
\left\{
\begin{array}{ll}
\va(r,t),\quad &\mbox{if }t\leq t_1,\\
\va_1(r,t-t_1)+\pi,\quad &\mbox{if }t> t_1.
\end{array}
\right.
\eeq
\endpf

\section{A crucial time to construct other solutions}
\setcounter{equation}{0}
\setcounter{theorem}{0}

By the argument in last section, we have constructed a global solution to the system \eqref{liquidcrystal}, which is regular except one point $(0,t_1)$. In this section, we would provide another time $t_2$ which is crucial in finding other solutions. To this end, we first construct a supersolution inspired by \cite{bph02} and \cite{topping02}. Define $\bar\lambda(t)>0$ satisfies the following equation
\beq\label{eqdel}
\bar\lambda'(t)=\delta e^{-2t}\bar\lambda^{\ve},\quad \bar\lambda(0)=0,
\eeq
for positive constants $\delta$ and $\ve\in(0,1)$. It is easy to see
$$
\bar\lambda^{1-\ve}(t)=\frac{\delta(1-\ve)}{2}\big(1-e^{-2t}\big)\geq0.
$$
Denote
$$\displaystyle
\tau(f):=f_{rr}+\frac{1}{r}f_r-\frac{\sin f\cos f}{r^2},
$$
and
\beq\label{nuq3.2}
\phi(r,\bar\lambda,t)=2\arctan \left(\frac{r}{e^t\bar\lambda}\right),
\eeq
which is a solution of $\tau(\phi)=0$ for any $t\in[0,+\infty)$ with $\phi(r,0,0)=\pi$ and $\phi(0,\bar\lambda,t)=0$ for $t>0$.
For any $\ve\in (0,1)$ and $\mu>0$, let $a=1+\ve$ and
\beq\label{nuq3.3}
\theta(r,\mu,t)=2\arctan \left(\frac{r^a}{e^{at}\mu}\right).
\eeq
Then $\theta(r,\mu,t)$ satisfies
\beq\label{nuq3.4}\displaystyle
\theta_{rr}+\frac{1}{r}\theta_r-\frac{a^2\sin\theta\cos\theta}{r^2}=0.
\eeq

\blm\label{lmsps}
There exist positive constants $\mu_0$ and $\delta_0$, such that for any $\mu\geq \mu_0$, $\ve\in(0,1)$ and $0<\delta\leq \delta_0$, the function defined by
\beq\label{nuq3.5}
\psi(r,t)=\phi(r,\bar\lambda(t),t)- \theta(r,\mu,t)
\eeq
satisfies
\begin{itemize}
\item[(i)]
$\psi(r,t)\in C^{1}((0,1)\times[0,+\infty))\cap C^{\infty}((0,1]\times(0,+\infty))\cap C^{1}([0,1]\times(0,+\infty))$.

\item[(ii)]
$\displaystyle \lim\limits_{r\rightarrow 0}\psi(r,0)=\pi,
\quad \lim\limits_{r\rightarrow 0}\psi(r,t)=0, \ t>0.$

\item[(iii)] $\psi$ is a supersolution of \eqref{drift-hm}, i.e.,
\beq\label{nuq3.6}
\displaystyle \psi_{rr}+\frac{\psi_r}{r}-\frac{\sin(2\psi)}{2r^2}-r\psi_r-\psi_t\leq 0.
\eeq
\end{itemize}
\elm

\pf
It is not so hard to verify the property (i) and (ii) by the definition of $\psi$. To prove (iii), choosing $\mu$ large enough so that $\theta(r,\mu,t)$ is small enough and
\beq\label{nuq3.7}
\cos \theta(r,\mu,t)\geq \frac{1}{1+\ve}
\eeq
for any $r\in[0,1]$ and $t\in[0,+\infty)$.
Direct caculation and using \eqref{nuq3.4} and \eqref{nuq3.7}, we have
\beq\label{nuq3.8}
\begin{split}
\tau(\psi)=&\frac{1}{r^2}\big[\sin\phi\cos\phi-\sin(\phi+\theta)\cos(\phi+\theta)-a^2\sin\theta\cos\theta\big]\\
=&\frac{1}{r^2}\big[\cos(2\phi-\theta)\sin\theta-a^2\sin\theta\cos\theta\big]\\
\leq&\frac{1}{r^2}\big[\cos(2\phi+\theta)\sin\theta-(1+\ve)\sin\theta\big]\\
\leq&-\frac{\ve\sin\theta}{r^2}
=-\frac{\ve}{r^2}\frac{2\mu e^{-at} r^a}{\mu^2+e^{-2at}r^{2a}}\\
\leq&-\frac{2\mu\ve e^{-at}}{\mu^2+1}r^{\ve-1}
\end{split}
\eeq
where we have also used the fact $\sin\theta\geq 0$ in first inequality and \eqref{nuq3.3} in second inequality.
From definitions \eqref{nuq3.2} and \eqref{nuq3.3}, we  have
\beq\label{nuq3.9}
\psi_t=-\frac{2re^t(\bar\lambda+\bar\lambda_t)}{e^{2t}\bar\lambda^2+r^2}+\frac{2a\mu r^{a}e^{at}}{\mu^2e^{2at}+r^{2a}}
\eeq
and
\beq\label{nuq3.10}
\psi_r=\frac{2\bar\lambda e^t}{e^{2t}\bar\lambda^2+r^2}-\frac{2a\mu r^{a-1}e^{at}}{\mu^2e^{2at}+r^{2a}}
\eeq
Combining (\ref{nuq3.9}) with (\ref{nuq3.10}) and using \eqref{eqdel}, we obtain
\beq\label{nuq3.11}
\psi_t+r\psi_r=-\frac{2re^t\bar\lambda_t}{e^{2t}\bar\lambda^2+r^2}=
-\frac{2\delta re^{-t}\bar\lambda^{\ve}}{e^{2t}\bar\lambda^2+r^2}.
\eeq
To conclude \eqref{nuq3.6}, it suffices to verify
\beq\label{nuq3.12}
\frac{2\delta re^{-t}\bar\lambda^{\ve}}{e^{2t}\bar\lambda^2+r^2}\leq \frac{2\gamma\mu e^{-at}\ve}{\mu^2+1}r^{\ve-1}.
\eeq
Let $s=\frac{r}{e^t\bar\lambda}$. Then \eqref{nuq3.12} is equivalent to
\beq\label{nuq3.13}
\frac{ s^{2-\ve}}{1+s^2}\leq \frac{\gamma\mu \ve}{\delta(\mu^2+1)},\ \forall s>0.
\eeq
It is easy to check that the function $\displaystyle\frac{ s^{2-\ve}}{1+s^2}$ has a maximum $M(\ve)$ depending only on $\ve$. Therefore, if we choose
$$\delta\leq \frac{\gamma\mu \ve}{M(\ve)(\mu^2+1)},$$
the inequality \eqref{nuq3.12} holds, which completes the proof of the lemma.
\endpf
\vspace{2mm}

Denote
$$
\overline\psi:=\psi+\pi
$$
where $\psi$ is defined by \eqref{nuq3.5} in Lemma \ref{lmsps} with fixed $\delta$ and large enough $\mu$ such that
$$
\pi-\arccos\left(\frac{\mu^2-1}{\mu^2+1}\right)\geq \beta,
$$
where $\beta\in(0,\pi)$ is given in \eqref{asspva0}. Now, we are ready for our main result on the crucial time $t_2$.

\begin{theorem}\label{thm4.1}
Let $\va_0$, $\tilde \va$ and $t_1$ be given as in Theorem \ref{thm3.2}. Then there exists $\ve>0$ and $t_1<t_2<+\infty$ such that
\beq\label{estt2}
\tilde\va(r,t)\leq \pi-\ve r,\quad \mbox{for } 0\leq r\leq 1,\ \ t\geq t_2.
\eeq
\end{theorem}

\pf
It is easy to check that $\overline\psi$ satisfies \eqref{nuq3.6}, and
$$
\overline\psi(r,0)\geq\tilde\va(r,t_1),\quad\mbox{for any }r\in[0,1].
$$
Combining the fact $\overline\psi(0,t)\geq\pi$ for all $t$ with the comparison principle, we conclude
$$
\overline\psi(r,t-t_1)\geq \tilde\va(r,t), \quad\mbox{for }t\geq t_1.
$$
Direct calculation indicates for any $r\in [0,1]$
$$
\frac{\partial}{\partial r}\overline\psi(r,t)=
\frac{2\bar\lambda e^t}{e^{2t}\bar\lambda^2+r^2}-\frac{2a\mu r^{a-1}e^{at}}{\mu^2e^{2at}+r^{2a}}
=\frac{2\bar\lambda e^{-t}}{\bar\lambda^2+e^{-2t}r^2}-\frac{2a\mu r^{a-1}e^{-at}}{\mu^2+e^{-2at}r^{2a}}
\rightarrow 0,\quad \mbox{as }t\rightarrow +\infty.
$$
Therefore, for any $\sigma>0$, there exists a time $t_\sigma>t_1$ such that
$$\tilde \va(r,t_\sigma)\leq \pi+\sigma r$$
for any $r\in[0,1]$.

In order to pick $t_2$, denote
$$
g(r,t)=\pi+l(\gamma-t) re^{-t}-lr^2e^{-2t},
$$
for $r\in[0,1]$, $t>0$, some constants $l>0$ and $1<\gamma<\ln 3$.
Direct calculation implies that $g(r,t)$ satisfies \eqref{nuq3.6} for $r\in[0,1]$ and $0<t<\ln 3$, while
\beq\label{nuq3.14}
g(1,0)>\pi,\quad g(1,t)>\beta,\ \ \mbox{for } 0<t<\ln 3,
\quad g_r(0,t)=l(\gamma-t)<0\ \ \mbox{for }\gamma<t<\ln 3,
\eeq
 with small enough $l>0$.
By the choice of $t_\sigma$, for small enough $\sigma>0$, there is $t_{\sigma}>0$ such that
$$
\tilde\va(r,t)\leq g(r,0), \quad \mbox{for }r\in[0,1].
$$
By the estimates \eqref{nuq3.14} and the comparison principle, it holds
$$
\tilde \va(r,t)\leq g(r,t-t_\sigma), \quad \mbox{for }t\in[t_\sigma,t_\sigma+\ln 3].
$$
In particular, setting $\ve=l$ and $\displaystyle t_2=t_\sigma+\frac{\gamma+\ln 3}{2}$, we conclude
$$
\tilde \va(r,t_2)\leq \pi-\ve r,\quad \mbox{for any }r\in[0,1].
$$
We can conclude \eqref{estt2} combining the fact that the function $g_1(r)=\pi-\ve r$ is also a supersolution, which completes the proof of theorem.
\endpf

\section{Proof of main theorems}
\setcounter{equation}{0}
\setcounter{theorem}{0}

This section is devoted to the construction of the second energy jump at some time beyond $t_2$ and the existence of another solution to complete the proof of Theorem \ref{mainth1}, and hence Theorem \ref{mainth1}.

\begin{theorem}\label{thm5.1}
Let $\va_0$, $\tilde \va$ and $t_1$ be given as in Theorem \ref{thm3.2}, and $t_2$ as in Theorem \ref{thm4.1}. For any $\tau\geq t_2$, there exist another solution $\va^{\tau}(r,t)$ to the system \eqref{drift-hm} in $(0,1)\times(0,+\infty)$, which is smooth except the points $(0,t_1)$ and $(0,\tau)$, and also satisfies
$$
\va^\tau(r,t)=\tilde\va(r,t),\quad \mbox{for any } (r,t)\in [0,1]\times[0,\tau],
$$
$$
\va^{\tau}(0,t)=\chi(t)=0,\quad\mbox{for } t>\tau.
$$
Furthermore, the energy $E(\va^{\tau}(t))$ is bounded and
\beq\label{engjp1}
\lim\limits_{t\downarrow \tau} E(\va^\tau(t))=4+E(\tilde \va(\tau)).
\eeq
\end{theorem}

\pf First consider the following initial and boundary value problem
\begin{equation}\label{nuq5.1}
\begin{cases}
\displaystyle\va^*_t+r\va^*_r=\va^*_{rr}+\frac{\va^*_r}{r}-\frac{\sin(2\va^*)}{2r^2}, \ 0< r< 1,\ t>0 \\
\varphi^*(r,0)=\tilde\varphi(r,\tau), \ 0\leq r\leq 1,\\
\varphi^*(0,t)=0, \ \ \varphi^*(1,t)=\beta, \ t>0,
\end{cases}
\end{equation}
and seek for a solution $\va^*\in C^{\infty}\big((0,1)\times[0,+\infty\big)\cap C\big([0,1]\times(0,+\infty\big)$. Notice that there is a singularity at $(0,0)$ since $\va^*(0,0)=\pi$ and $\va^*(0,t)=0$ for $t>0$. Thus, we need first construct an approximate sequence of initial data $\va_0^n(r)$ for any $n\in\mathbb N$, i.e., for any $r_n\in (0,1)\rightarrow 0$ as $n\rightarrow+\infty$, define a family of functions
\beq\notag
\va_0^n(r):=\left\{
\begin{array}{ll}
\displaystyle 2\arctan \left(\frac{r}{r_n}\tan\left(\frac{\tilde \va(r_n,\tau)}{2}\right)\right),\quad &\mbox{if }0\leq r\leq r_n,\\
\displaystyle \tilde \va(r,\tau),\quad &\mbox{if }r_n<r\leq 1.
\end{array}
\right.
\eeq
By the estimate \eqref{estt2}, $\tilde\va(r_n,\tau)<\pi$ and
$\tilde\va(r_n,\tau)\uparrow \pi$ as $t_n\rightarrow 0$, so the functions $\va_0^n(r)$ are well-defined for any $n\in\mathbb N$ and smooth in $[0,1]$ except the point $r_n$. It is direct to check that
$\va_0^n(0)=0$ and $\va_0^n(r)\leq \tilde \va(r,\tau)$ in $[0,1]$.
By the results in \cite{vdh01}, it holds
$$
E_{0,r_n}=2(1-\cos (\tilde \va(r_n,\tau))).
$$
Similar to the proof of \eqref{engdrop1}, we conclude
\beq\label{nuq5.2}
\lim\limits_{n\rightarrow+\infty}E(\va_0^n)
=4+E(\tilde\va(\tau)).
\eeq
By the standard diagonal procedure of approximation of a function in $H^1$ by $C^{\infty}$ function, we may assume $\va_0^n\in C^{\infty}[0,1]$, which still satisfies all the properties above and in addition, $\va_0^n\rightarrow \tilde\va(\tau)$ uniformly in $C^k$ in any compact subsets of $(0,1]$ for any $k\in\mathbb N$.

Let $\psi(r,t)$ be the function defined as \eqref{nuq3.5}. By the Theorem \ref{thm4.1} and Lemma \ref{lmsps}, choose suitable $\mu>0$ and $\delta>0$ such that
\beq\label{nuq5.3}
\tilde\va(r,\tau)\leq \psi(r,0),\quad \mbox{for }r\in[0,1],
\eeq
and
\beq\notag
\psi(1,t)\geq \beta,\quad\mbox{for }t\in[0,t'],
\eeq
with some small $t'>0$. By the local existence result in \cite{hllw16}, the system \eqref{nuq5.1} with initial data $\va_0^n$ has a unique smooth solution $\va^n$ in $[0,1]\times[0,t'')$, where $t''\geq t'$ without loss of generalization. Combining the fact $0\leq \va_0^n(r)\leq \tilde \va(r,\tau)$ with the comparison principle Lemma \ref{complm}, for any $n\in\mathbb N$ it holds
\beq \label{nuq5.4}
0\leq \va^n(r,t)\leq \psi(r,t),\quad \mbox{for any }(r,t)\in [0,1]\times[0,t'].
\eeq
Choose small enough $\mu^*$ such that
$$
\psi^*:=2\arctan\left(\frac{r}{\mu^*}\right)\geq\psi(r,t'),\quad\mbox{for }r\in[0,1].
$$
Direct calculation implies
\beq\notag
\displaystyle \psi^*_{rr}+\frac{\psi^*_r}{r}-\frac{\sin(2\psi^*)}{2r^2}-r\psi^*_r-\psi^*_t\leq 0.
\eeq
By all these properties and comparison principle, we have
\beq\label{nuq5.5}
0\leq \va^n(r,t)\leq \overline\psi(r,t):=
\left\{
\begin{array}{ll}
\psi(r,t),\quad&\mbox{if }0\leq t\leq t',\\
\psi^*(r,t)\quad&\mbox{if }t'<t\leq t''.
\end{array}
\right.
\eeq
By similar argument as in \cite{hllw16} Section 4, the solution $\va^n$ can be extended beyond $t''$ smoothly, which implies $t''=+\infty$.

Now, we proceed to take the limit as $n\rightarrow+\infty$. For any $0<T<+\infty$ denote
$$
W_T:=\big\{f\big|~ f\in L^{\infty}((0,T); H^1_r(0,1))\cap L^{\infty}([0,1]\times[0,T)),\  f_t\in L^2((0,T);L^2_r(0,1))\big\}
$$
By the energy estimates in Remark \ref{rmk2.3}, $\va^n$ is uniform bounded in $W_T$. There exists a function $\va^*\in W$ such that
$$
\va^*(r,t)\leq \overline \psi(r,t),\quad \mbox{for any }(r,t)\in[0,1]\times[0,T],
$$
and choosing subsequences if necessary, as $n\rightarrow +\infty$,
$$
\va^n\rightarrow \va^*,\quad\mbox{weakly in } W_T,
$$
$$
\va^n\rightarrow \va^*,\quad\mbox{strongly in } L^2([0,1]\times[0,T]).
$$
And $\va^n\rightarrow \va^*$ uniformly on any compact subsets of $(0,1]\times[0,T]$ or $[0,1]\times(0,T]$. It is not hard to see  $\va^*\in C^{\infty}\big((0,1)\times[0,T]\big)\cap C\big([0,1]\times(0,T]\big)$ is a solution to \eqref{nuq5.1}. By the energy estimate in Remark \ref{rmk2.3}, \eqref{nuq5.2} and the lower-semicontinuity of energy functional $E(\cdot)$, we conclude $E(\va^*(t))$ is bounded for any $t>0$.

Now, we are ready to prove the energy jump \eqref{engjp1}. This can be done similarly as \eqref{engdrop1}. For a small enough $r'\in(0,1)$
\beq\notag
\begin{split}
&\lim\limits_{t\downarrow 0} E(\va^*(t))\\
=&\lim\limits_{t\downarrow 0}\big(E_{0, r'}(\va^*(t))+E_{r',1}(\va^*(t))\big)\\
\geq& \lim\limits_{t\downarrow 0}\big(2(1-\cos \va^*(r', t))+E_{r',1}(\va^*(t))\big)\\
=&2(1-\cos \tilde\va(r',\tau))+E_{ r',1}(\tilde\va(\tau)).
\end{split}
\eeq
Taking the limit $r'\rightarrow 0$ implies \eqref{engjp1}.
Therefore, $\tilde\va(r,t)$ can be defined as follows
\beq\notag
\va^\tau(r,t):=
\left\{
\begin{array}{ll}
\tilde\va(r,t),\quad &\mbox{if }t\leq \tau,\\
\va^*(r,t-\tau)+\pi,\quad &\mbox{if }t> \tau,
\end{array}
\right.
\eeq
which completes the proof of Theorem \ref{thm5.1}.
\endpf

\bigskip
%

\end{document}